%
\magnification=\magstep1
\input amstex
\UseAMSsymbols
\input pictex
\NoBlackBoxes
   
   \font\rmk=cmr8    \font\itk=cmti8  \font\ttk=cmtt8

   \newcount\notenumber
   
   \def\note{\advance\notenumber by 1 
       \plainfootnote{$^{\the\notenumber}$}}  
\def\ab{\operatorname{ab}}

\def\Im{\operatorname{Im}}

\def\cok{\operatorname{cok}}
\def\stat{\operatorname{stat}}
\def\adstat{\operatorname{adstat}}
\def\mod{\operatorname{mod}}
\def\Mod{\operatorname{Mod}}

\def\Hom{\operatorname{Hom}}

\def\End{\operatorname{End}}
\def\Ext{\operatorname{Ext}}

\def\Ker{\operatorname{Ker}}

\def\add{\operatorname{add}}

\def\bdim{\operatorname{\bold{dim}}}
\def\op{{\text{op}}}

\def\arr#1#2{\arrow <1.5mm> [0.25,0.75] from #1 to #2}

\vglue1truecm
\centerline{\bf Static subcategories}
\medskip
\centerline{\bf  of the module category of a finite-dimensional hereditary algebra.}
		\bigskip\medskip
\centerline{Mustafa A\. A\. Obaid, S\. Khalid Nauman,}
		    \smallskip
\centerline{Wafaa M\. Fakieh and Claus Michael Ringel (Jeddah)}
		  \bigskip\bigskip

\plainfootnote{}
{\rmk 2010 \itk Mathematics Subject Classification. \rmk 
Primary:
        16G20,  
        16G60, 
        05E10. 
Secondary:
        16D90, 
        16G70. 
\newline Key words and phrases: Finite-dimensional hereditary algebras. Representation types: 
tame and wild, strictly wild. 
$M$-static modules, $M$-adstatic modules. The smallest exact abelian subcategory containing a given module. Bricks. $\ab$-projective modules. Nakayama algebras. Triple modules. 
}

{\narrower\narrower Abstract. Let $k$ be a field,
$\Lambda$ a finite-dimensional hereditary $k$-algebra and $\mod\Lambda$ the category
of all finite-dimensional $\Lambda$-modules. 
We are going to characterize the representation type of $\Lambda$ (tame or wild)
in terms of the possible subcategories $\stat M$ of all $M$-static modules,
where $M$ is an indecomposable $\Lambda$-module. 
\par}
	\bigskip\bigskip

{\bf 1. Introduction.}
	\medskip 
Let $k$ be a field and $\Lambda$ a finite-dimensional $k$-algebra. 
The $\Lambda$-modules to be considered will be assumed to be finite-dimensional left modules.
Let $\mod\Lambda$ be the category of these $\Lambda$-modules. 

Let $M$ be a $\Lambda$-module, 
let $\Gamma(M) = \End(M)^\op$
be the opposite of the endomorphism ring of $M$. 
Consider the functor
$F= \Hom_\Lambda(M,-):\mod\Lambda \to \mod\Gamma(M)$ and its left adjoint $G = M\otimes_{\Gamma(M)}-$
(usually, we will just write $\Hom(M,-)$ instead of $\Hom_\Lambda(M,-)$ and $M\otimes-$ instead of
$M\otimes_{\Gamma(M)}-$). 
There are canonical maps
$$
\align
 \mu_N: GF(N) \to N,\quad & \text{for $N$ a $\Lambda$-module, and} \cr
 \nu_X: X \to FG(X),\quad & \text{for $X$ a $\Gamma(M)$-module.} \cr
\endalign
$$ 
A $\Lambda$-module $N$ is said to be {\it $M$-static} provided $\mu_N$ is an isomorphism
and we denote by $\stat M$ the subcategory of $\mod\Lambda$ given by all
$M$-static modules (subcategories considered in the paper are assumed to be full).
A $\Gamma(M)$-module $X$ 
is said to be {\it $M$-adstatic} provided $\nu_X$ is an isomorphism
and we denote by $\adstat M$ the subcategory of $\mod\Gamma(M)$ given by all
$M$-adstatic modules. 
It is easy to see that the functors $F$ and $G$ provide
an equivalence between the categories $\stat M$ and $\adstat M$. 
For general properties of static and adstatic modules we may refer to Alperin [A], Nauman [N1], [N2] and 
Wisbauer [W].

We denote by $\add M$ the subcategory of $\mod\Lambda$
given by all direct summands of (finite) direct sums of copies of $M$. 
Let $\cok(M)$ be the subcategory of all cokernels of maps in $\add M$, one always
has
$$
 \add M \subseteq \stat M \subseteq \cok(M).
$$

Recall that a ring is said to be {\it hereditary} provided submodules of projective modules are
projective. In this paper, we will deal with finite-dimensional  hereditary $k$-algebras
and we want to characterize the representation type of such an algebra $\Lambda$ by looking at
the static subcategories of $\mod\Lambda$. These algebras have been studied thoroughly (see
in particular [DR1] and [R]; we will recall the relevant facts in section 7): 
Such an algebra is either tame or wild. In case $\Lambda$ is tame, the category
$\mod\Lambda$ consists of directed components and separating tubes.
In case $\Lambda$ is wild, one knows that $\Lambda$ is strictly wild.
We will use this knowledge in order to 
characterize the tameness of $\Lambda$ in terms of the subcategories $\stat M$ with $M$ indecomposable. 
	\bigskip
Given a module $M$, let $\ab M$ be the smallest exact 
abelian subcategory in $\mod\Lambda$ which contains $M$,
it is the intersection of all exact abelian subcategories containing $M$, and also
the closure of $\add M$ using (inductively) kernels and 
cokernels. If $M$ is an indecomposable module, then $\ab M = \add M$ if and only if
$M$ is a brick 
(a {\it brick} is a module whose endomorphism ring is a division ring),

We say that $M$ is {\it $\ab$-projective} provided $M$, considered as an object of
the abelian category $\ab M$, is projective. Here are typical examples: of course, 
all bricks are $\ab$-projective; second, if $I$ is an ideal of $\Lambda$
which annihilates $M$ and $M$ is projective as a $\Lambda/I$-module, then
$M$ is $\ab$-projective.
	\bigskip
Recall that a module is said to be {\it serial} provided it has a unique composition series. 
A finite dimensional $k$-algebra is said to be a
{\it Nakayama algebra} provided any indecomposable module is serial 
(it is sufficient to assume that any indecomposable projective
and any indecomposable injective module is serial). 
	\bigskip 
{\bf Theorem 1.} {\it Let $\Lambda$ be a finite-dimensional hereditary $k$-algebra.
The following conditions are equivalent:
\item{\rm (i)} $\Lambda$ is tame.
\item{\rm(ii)} Any indecomposable module is $\ab$-projective.
\item{\rm(iii)} If $M$ is indecomposable, then $\Gamma(M) = \End(M)^\op$ is a Nakayama algebra 
  \newline and 
    $\adstat M = \mod\Gamma(M).$     
\item{\rm(iv)} If $M$ is indecomposable, then $\stat M$ is abelian.
\item{\rm(v)} If $M$ is indecomposable, then $\stat M = \cok(M).$}
	\bigskip 
{\bf Theorem 2.} {\it Let $\Lambda$ be a finite-dimensional hereditary $k$-algebra. 
The following conditions are equivalent:
\item{\rm (i)} $\Lambda$ is wild.
\item{\rm(ii)} There exists an indecomposable module $M$ which is not a brick, but 
    $\stat M = \add M.$
\item{\rm(iii)} There exists a finite extension field $k'$ of $k$ such that
     for any finite-dimensional $k'$-algebra $\Gamma$, there is a module $M$ such that
     $\stat M$ is equivalent to $\mod\Gamma.$\par}

	\bigskip\bigskip 
{\bf 2. $M$-static modules.}
	\medskip
We will use a well-known characterization (see [A] and [W])
of the $M$-static modules as cokernels of maps in $\add M$.
We need the following definitions: Let $M, N$ be $\Lambda$-modules. A map $q:M' \to N$
is called a {\it right $M$-approximation} provided $M'$ belongs to $\add M$
and for any map $g:M \to N$ there is $g':M \to M'$ such that $g = qg'$ (this just means that
$\Hom(M,q)$ is surjective). 
A minimal right $M$-approximation is a right $M$-approximation which is right
minimal. In case $\Lambda$ is a finite-dimensional algebra and $M,N$ are
finite-dimensional
$\Lambda$-module, a minimal right $M$-approximation of $N$ exists and
we denote by $\Omega_M(N)$ its kernel (see [DR2]), it is unique up to isomorphism.  
	\medskip
Here is the characterization in the case of dealing with finite-dimensional modules for
finite-dimensional algebras:
	\medskip 
{\bf Proposition 1.} {\it Let $\Lambda$ be a finite-dimensional algebra and $M,N$
finite-dimensional $\Lambda$-modules. The following assertions are equivalent:
\item{\rm (i)} $N$ is $M$-static. 
\item{\rm (ii)} Both $N$ and $\Omega_M(N)$ are generated by $M$.
\item{\rm (iii)} There is an exact sequence
$$
 M'' @>>> M' @>q>> N @>>> 0
$$ 
with $M',M''$ in $\add M$, such that $q$ is a right $M$-approximation.
\item{\rm (iv)} There is an exact sequence
$$
 M'' @>>> M' @>>> N @>>> 0
$$ 
with $M',M''$ in $\add M$, such that the sequence remains exact when we apply $\Hom(M,-).$
\par}
	\medskip
{\bf Remark.} Let us stress that the sequences provided in condition 
(iii) may not remain exact when we apply $\Hom(M,-)$ as required in (iv).
Here is an example: Take the quiver with two vertices $1,2$, two arrows $\alpha,\beta\:1 \to 2$ and
one arrow $\gamma\:2 \to 1$, and take as relations the paths $\alpha\gamma, \beta\gamma$ 
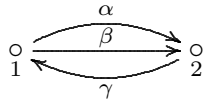
$$
\hbox{\beginpicture
	\setcoordinatesystem units <2cm,1cm>
\put{\beginpicture
\multiput{$\circ$} at 0 0  1 0 /
\put{$\ssize 1$} at 0 -.2
\put{$\ssize 2$} at 1 -.2
\put{$\ssize \alpha$} at .5 .45
\put{$\ssize \beta$} at .5 .15
\put{$\ssize \gamma$} at .5 -.45
\setquadratic
\plot 0.9 0.1 0.5 0.3  0.1 0.1 /
\plot 0.9 -.1  0.5 -.3  0.1 -.1 /
\setlinear
\plot 0.9 0  0.1 0 /
\arr{0.9 0.11}{0.91 0.1}
\arr{0.9 0}{0.91 0}
\arr{0.1 -.11}{0.09 -.1}
\endpicture} at 0 0
\endpicture}
$$
We consider the module $M = I(1)$ (the indecomposable injective module with socle $S = S(1)$). 
Let $N = M/S$ (this is the 3-dimensional indecomposable injective module for the
Kronecker quiver with arrows $\alpha,\beta$). Let $f\:M \to M$ be an endomorphism of $M$ with image $S$, and
$q\:M\to N$ its cokernel, thus we deal with the exact sequence
$$
 M @>f>> M @>q>> N @>>> 0. \tag$*$
$$ 
Since $\dim\Hom(M,N) = 1$, we see that $q$ is a right $M$-approximation, thus the sequence is
as required in (iii). However, $\dim \Hom(M,S) = 2$ and $\Hom(M,S)$ is annihilated by the radical of
$\Gamma(M)$, thus a minimal right $M$-approximation of $S$ is of the form $M^2 \to S.$
If we write $f = up$ with $p\:M\to S$ and choose $p'\:M \to S$ such that $p,p'$ is a basis
of $\Hom(M,S)$, then $up'\in \Hom(M,M)$ is in the  kernel of the map $\Hom(M,q)$, but not
in the image of $\Hom(M,f).$ This shows that $(*)$ does not stay exact when we apply $\Hom(M,-)$.
	\medskip

Here is the proof of Proposition 1 (see [A], [W]). 

(i) $\implies$ (ii). Since $\mu_N:M\otimes\Hom(M,N) \to N$ is surjective, 
we easily see that $N$ is generated by $M$. Thus, the minimal right $M$-approximation
$q\:M'\to N$ is surjective and we have an exact sequence
$$
 0 @>>> \Omega_M(N) @>u>> M' @>q>> N @>>> 0.
$$
The functor $\Hom(M,-)$ is left exact, thus 
$$
 0 @>>> \Hom(M,\Omega_M(N)) @>>> \Hom(M,M') @>\Hom(M,q)>> \Hom(M,N)
$$
is exact. Since $q$ is a right $M$-approximation, the map $\Hom(M,q)$ is surjective,
thus we deal with the exact sequence
$$
 0 @>>> \Hom(M,\Omega_M(N)) @>>> \Hom(M,M') @>\Hom(M,q)>> \Hom(M,N) @>>> 0.
$$
The functor $M\otimes -$ is right exact, thus the 
upper sequence in the following commutative diagram is exact (as is the lower one):
$$
\CD
 @. \ssize M\otimes \Hom(M,\Omega_M(N)) @>>> \ssize M\otimes \Hom(M,M') @>\sssize M\otimes \Hom(M,q)>>\ssize  M\otimes \Hom(M,N) @>>> 
   \ssize  0 \cr
 @.     @V\mu_{\Omega_M(N)}VV      @V\mu_{M'}VV                   @V\mu_N VV  \cr
 \ssize 0 @>>> \ssize \Omega_M(N) @>u>>   \ssize        M'    @>q>>             \ssize   N    @>>> \ssize  0
\endCD
$$
The maps $\mu_{M'}$ and $\mu_N$ are isomorphisms, thus it follows that
the map $\mu_{\Omega_M(N)}$ is surjective. But this means that $\Omega_M(N)$ is generated by $M$.

(ii) $\implies$ (iv). Since $N$ is generated by $M$, there is an exact sequence
$$
 0 @>>> \Omega_M(N) @>u>> M' @>q>> N @>>> 0
$$
where $q$ is a minimal right $M$-approximation of $N$. Since $\Omega_M(N)$ 
is generated by $M$, and right $M$-approximation $p\:M'' \to \Omega_M(N)$
is surjective. Thus we obtain an exact sequence
$$
  M'' @>up>> M' @>q>> N @>>> 0
$$
and since $u$ is injective and $p,q$ are right $M$-approximations, it follows that
this sequence is mapped under $\Hom(M,-)$ to an exact sequence. 

(iv) $\implies$ (i). We start with the exact sequence
$$
 M'' @>>> M' @>>> N @>>> 0
$$ 
and apply $\Hom(M,-)$. By assumption, we obtain the exact sequence
$$
 \Hom(M,M'') @>>> \Hom(M,M') @>\Hom(M,q)>> \Hom(M,N) @>>> 0.
$$
We apply the right exact functor $M\otimes-$ and obtain the upper exact sequence of the
following commutative diagram:
$$
\CD
  M\otimes\Hom(M,M'') @>>> M\otimes\Hom(M,M') @>M\otimes\Hom(M,q)>> M\otimes\Hom(M,N) @>>> 0 \cr
   @V\mu_{M''}VV      @V\mu_{M'}VV                   @V\mu_N VV  \cr
  M'' @>u>>          M'    @>q>>               N    @>>> 0
\endCD
$$
Since the vertical maps $\mu_{M''}$ and $\mu_{M'}$ are isomorphisms, also
$\mu_n$ is an isomorphism.
	\medskip 
The implications (iv) $\implies$ (iii) and (iii) $\implies$ (ii) are trivial.  

\hfill $\square$
	\bigskip
{\bf Proposition 2.} {\it Let $M$ be an indecomposable  module such that $\Gamma(M)$ is
a Nakayama algebra. If $N$ is an indecomposable $M$-static $\Lambda$-module, then there is a submodule
$U$ of $M$ such that $N$ is isomorphic to $M/U$. If $U$ is a submodule of $M$ and $p\:M \to M/U$ is
the canonical projection, then $M/U$ is $M$-static if and only if $p$ is a right $M$-approximation and
$U$ is the image of an endomorphism of $M$.}
	\medskip
Proof. First, assume that $U$ is a submodule of $M$ such that the canonical projection $M \to M/U$ is a right
$M$ approximation. Then the implication (iii) $\implies$ (i) in Proposition 1 shows that $M/U$
is $M$-static. 

Conversely, let us assume that $N$ is indecomposable and $M$-static. We have to show that 
$N$ is of the form $M/U$ where $U$ is a submodule of $M$, and that for any submodule $U$ of $M$ such
that $M/U$ is isomorphic to $N$, the projection map $M \to M/U$ is a right
$M$-approximation and that $U$ is the image of an endomorphism of $M$. 
These assertions are trivially true in case $N = M$, 
thus we will assume that $N$ is not isomorphic to $M$. According to Proposition 1, 
there exists an exact sequence
$$
 M'' @>f>> M' @>q'>> N @>>> 0 \tag{$*$} 
$$ 
with $M', M''$ in $\add M$, such that $q$ is a right $M$-approximation. Since $M$ is indecomposable,
$M' = M^a, M'' = M^b$ for some natural numbers $a,b$. The map $f\:M^b \to M^a$ is given by an
$(a\times b)$-matrix $C$ with coefficients in $\Gamma(M)$. Since $\Gamma(M)$ is a Nakayama algebra,
there are invertible square matrices $A, B$ such that $ACB$ is a diagonal matrix (in order to see this,
one can use the usual matrix reduction as in the case of matrices with coefficients in a field;
of course, the assertion corresponds to the fact that all $\Gamma(M)$-modules are direct sums of
serial modules: the matrix $C$ describes a projective presentation of a $\Gamma(M)$-module).
We use the matrices $A$ and $B$ in order to define automorphisms of $M^a$ and $M^b$. Thus, without
loss of generality we can assume that $f\:M^b \to M^a$ is given by a diagonal matrix. One of the
diagonal coefficients of this matrix, say the coefficient $c$ at the position $(1,1)$ has to be non-invertible. 
Since $N$ is indecomposable, it follows that $q'$ vanishes on $0\oplus M^{a-1}$ and we denote
by $q$ the restriction of $q'$ to $M\oplus 0$. Since $q'$ is a right $M$-approximation of $N$,
also $q$ is a right $M$-approximation of $N$. Since $N$ is not isomorphic to $M$, we must have
$c \neq 0$. It follows that the sequence $(*)$ splits off 
the exact sequence
$$
 M @>c>> M @>q>> N @>>> 0. \tag{$**$}
$$ 
This shows that $N$ is isomorphic to $M/U$, where $U$ is the image of $c$.
Let $U$ be a submodule of $M$ with an isomorphism $g\:M/U \to N$. Let $p\:M \to M/U$ be
the canonical projection. Since $q$ is a right $M$-approximation, there is a map
$g\:M \to M$ such that $qg' = gp$. Since $\Gamma(M)$ is a local Nakayama algebra, and $c,g'$ are
endomorphisms of $M$, we have either $g'(M) \subseteq c(M)$ or $c(M) \subseteq g'(M)$.
Now $g'(M) \subseteq c(M)$ is impossible, since $qc = 0$, whereas  the image of $qg' = gp$
is equal to $N$, thus non-zero. This shows that $c(M) \subseteq g'(M).$ Since $qg'$ is surjective,
we see that also $g'$ is surjective, thus $g'$ is an isomorphism. As a consequence, the pair $(g',g)$
is an isomorphism from $p\:M \to M/U$ to $q\:M \to N$. This shows that $p$ is a right
$M$-approximation of $M/U$. Also, $g'$ maps the kernel $U$ of $p$ onto the kernel of $q$ and
this is the image of the endomorphism $c$, thus $U$ is the image of an
endomorphism of $M$. \hfill $\square$
	\bigskip\bigskip 

{\bf 3. The subcategories $\ab M$.}
	\medskip
Looking at a module $M$, the subcategory $\ab M$ should be seen as an important invariant.
But it seems that a study of this invariant has been neglected up to now, thus we want
to provide at least some basic properties. In order to do so, let us first consider the general
setting of dealing with an arbitrary ring $R$ and any
(left) $R$-module $M$. Given a ring
$R$, let $\Mod R$ be the category of all $R$-modules. 

If $M$ is an $R$-module, we define $\ab M$ as the smallest exact abelian
subcategory\plainfootnote{$^1$}{Remark:
We should stress that $\ab M$ refers to ``smallest exact abelian subcategory containing $M$'', 
not to ``smallest abelian subcategory containing $M$'': 
note that the latter formulation would not even make sense,
since {\it the intersection of abelian subcategories which contain $M$
is not necessarily abelian.} Here is an example: Consider the quiver $1 \leftarrow 2 \leftarrow 3$,
and let $M_1 = S(1), M_2 = I(1), M_3 = I(2)$ and $M_3' = S(3)$ (for any vertex $x$ of a quiver without loops, 
we denote by $S(x)$ the simple module concentrated at $x$, and by $I(x)$ the injective envelope of $S(x)$).
Let $M = M_1\oplus M_2$.
Then both $M_1,M_2,M_3$ and $M_1,M_2,M_3'$ are the
indecomposable objects in abelian subcategories $\Cal A$ and $\Cal A'$.
Both $\Cal A$ and $\Cal A'$ contain $M$, and the intersection of $\Cal A$ and $\Cal A'$ is $\add M,$ 
but $\add M$ is not abelian.} 
of the category $\Mod R$ of all $R$-modules which contains $M$, it is the intersection
of all exact abelian subcategories containing $M$. By definition, a subcategory $\Cal C$ of an abelian
category $\Cal A$ is an {\it exact abelian subcategory} provided $\Cal C$ is closed under kernels, 
cokernels and
direct sums (and this is equivalent to say that $\Cal C$ itself is an abelian category and 
the inclusion functor $\Cal C \to \Cal A$ is exact).

 Let us show that {\it $\ab M$ is the closure of $\add M$
using kernels and cokernels.} If $\Cal C$ is a  subcategory
of $\Mod R$, define inductively subcategories $\ab_n(\Cal C)$ of $\Mod R$ as follows: 
Let $\ab_0M = \add M.$ If $\ab_{n-1}M$ is already defined for some $n\ge 1$, then let 
$\ab_nM$ be the subcategory of $\Mod R$ given by all $R$-modules which are kernels or cokernels
of maps in $\ab_{n-1}M$. Note that $\ab_{n-1}M \subseteq \ab_nM$ and that $\ab_nM$ is closed
under direct sums. We obtain an increasing chain of subcategories closed under
direct sums
$$
 \add M = \ab_0M \subseteq \cdots \subseteq \ab_{n-1}M \subseteq \ab_nM \subseteq \cdots .
$$
We denote by $\ab_\omega M = \bigcup_n\ab_nM$ the union, this is a subcategory of $\Mod R$ which is closed under
kernels, cokernels and direct sums, thus an exact abelian subcategory of $\Mod R.$
Of course, if $\Cal A$ is any exact abelian subcategory of $\Mod R$
which contains $M$, then $\ab_\omega M \subseteq \Cal A.$ This shows that $\ab M = \ab_\omega M.$
	\bigskip 
A {\it length category} is by definition an abelian category $\Cal A$ such that every object in
$\Cal A$ has a finite composition series (thus finite length). Given an object $A$ in a length category
$\Cal A$, the {\it Loewy length} of $A$ is the smallest number $t\ge 0$ such that $A$ has a filtration
with $t$ semisimple factors. 
	\medskip 
{\bf Proposition 3.} {\it Let $R$ be any ring. 
If $M$ is an $R$-module of finite length $t$, then
$\ab M$ is a length category with at most $t$ simple objects such that all objects in
$\ab M$ have Loewy length at most $t$.

If $\Lambda$ is a finite-dimensional $k$-algebra and $M$ a finite-dimensional $\Lambda$-module, 
then there is 
a finite-dimesnional $k$-algebra $\Lambda'$ such that $\ab M$ is equivalent to $\mod\Lambda'$.}
	\bigskip
Proof: First, let $R$ be any ring and $M$ an $R$-module of finite length $t$.
Consider $M$ as an object in the abelian subcategory $\ab M$. Since any subobject of $M$ 
is an  $R$-submodule,
we see that $M$, considered as an object of $\ab M$, has (relative) length at most $t$, thus also
(relative) Loewy length at most $t$. 
Let $S_1,\dots, S_s$ be the (relative) simple objects of $\ab M$ which occur as 
factors in a (relative) composition series of $M$,
then $s\le t$. Note that the objects $S_1,\dots, S_s$ are pairwise orthogonal bricks in $\mod\Lambda$,
the process of simplification (see [R]) shows that the subcategory $\Cal F(S_1,\dots,S_s)$
of all $\Lambda$-modules with a filtration with factors of the form $S_1,\dots, S_s$ is
an exact abelian subcategory of $\Mod R$ whose (relative) simple objects are precisely
the objects $S_1,\dots,S_s$. We also may consider the subcategory $\Cal F(S_1,\dots,S_s;t)$
of $R$-modules $N$ with a filtration with factors $S_1,\dots, S_s$, such that the (relative)
Loewy length of $N$ is at most $t$.
Then $\Cal F(S_1,\dots,S_s;t)$
is an exact abelian subcategory which contains $M$, thus $\ab M \subseteq 
\Cal F(S_1,\dots,S_s;t)$. Since the modules $S_1,\dots,S_s$ belong to $\ab M$, we see
that these are precisely the (relative) simple objects of $\ab M$ and that any object of $\ab M$
has (relative) Loewy length at most $t$. 

Now assume that $\Lambda$ is a finite-dimensional $k$-algebra, where $k$ is 
a field. Let $M$ be a $\Lambda$-module of finite length. The category $\mod\Lambda$
is both $\Hom$-finite and $\Ext$-finite. 
If $N,N'$ belong to $\ab M$, then $\Ext^1_{\ab M}(N',N)$ is a subspace of $\Ext^1_\Lambda(N',N)$,
thus finite-dimensional. Thus $\ab M$ is a length $k$-category which is $\Hom$-finite and $\Ext$-finite,
with finitely many simple objects and bounded Loewy length. It is well-known that such a
category has a progenerator, say $P$. If $\Lambda' = \End(P)^\op$, then $\ab M$ is 
equivalent to $\mod\Lambda'$.\par
 \hfill $\square$
	\bigskip
Examples to have in mind: 
	\medskip
(1) If $M$ is a brick, then $\ab M = \add M$.
More generally, if $M$ is the direct sum of pairwise orthogonal bricks, then $\ab M = \add M$.
	\medskip 
(2) Let $\Lambda$ be the Kronecker algebra, this is the path algebra of the quiver with two
vertices $1,2$ and two arrows $1\to 2$, the $\Lambda$-modules are usually called {\it Kronecker modules.}
A Kronecker module is said to be regular provided it is the direct sum of indecomposable modules
with even dimension. The regular Kronecker modules
form an exact abelian subcategory $\Cal R$ of $\mod\Lambda$ which has infinitely many 
(relative) simple objects. Also, note that the (relative) Loewy length of the objects in $\Cal R$ 
is not bounded. 

Let $M$ be an indecomposable Kronecker module. If $M$ is not regular, then $M$ belongs to
the preprojective or the preinjective component, and this implies that $M$ is a brick, thus
$\ab M = \add M.$ Thus, assume that $M$ is regular. Then there exists a (relative) simple regular 
object $X$ such that $M$ has a filtration with all factors isomorphic to $X$. Assume that this
filtration (it is unique) has length $e$. Then the indecomposable regular modules with a filtration
with at most $e$ factors of the form $X$ form an exact abelian subcategory, and this is just
$\ab M.$
	\medskip
(3) Assume that $\Lambda$ is a finite-dimensional algebra and that $M$ is a faithful $\Lambda$ module.
If $\ab M$ contains all simple $\Lambda$-modules, then $\ab M = \mod\Lambda$. Namely, if $\ab M$ contains
all simple $\Lambda$-modules, then clearly $\ab\Lambda$ is closed under submodules. On the other hand,
if $M$ is faithful, then 
there is an embedding ${}_\Lambda\Lambda \subseteq M^t$ for some natural number $t$. Thus, it follows that
${}_\Lambda\Lambda$ is contained in $\ab M$. Since any $\Lambda$-module $N$ has a free presentation,
we see that $N$ belongs to $\ab M.$
	\medskip 
(4) Whereas for a brick $M$, the abelian category $\ab M$ has a unique simple object
(and no other indecomposable objects), already for $\dim\End(M) = 2,$ the category $\ab M$ may have
arbitrarily many simple objects. Here is an example: We consider the algebra $\Lambda$ given by the
following quiver
$$
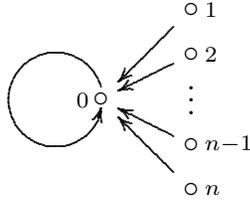
\hbox{\beginpicture
	\setcoordinatesystem units <1cm,1cm>
\multiput{} at 0 -.5  2 1 /
\multiput{$\circ$} at  0 0  1 1   1 0.5   1 -0.5  1  -1 /
\arr{0.8 0.8}{0.2 0.2}
\arr{0.8 0.4}{0.2 0.1}
\arr{0.8 -.8}{0.2 -.2}
\arr{0.8 -.4}{0.2 -.1}
\circulararc 330 degrees from 0 0.15 center at -.5 0 
\arr{-.03 -.2}{0 -.1}
\put{$\ssize 1$} [l] at 1.15 1 
\put{$\ssize 2$} [l] at 1.15 0.5 
\put{$\ssize n-1$} [l] at 1.15 -.5 
\put{$\ssize n$} [l] at 1.15 -1
\put{$\ssize 0$} at -.2 0 
\put{$\vdots$} at 1 0.1 
\endpicture}
$$
with all paths of length 2 as relations. Let $M = I(0).$ Then $\dim\End(M) = 2$.
If $f$ is a non-zero nilpotent endomorphism of $M$, then the cokernel of $f$ is the direct sum
of all simple modules. Thus $\ab M$ contains all simple $\Lambda$-modules. (Since $M$ is faithful,
it follows that $\ab M = \mod\Lambda$.) 
	\bigskip\bigskip
{\bf 4. $\ab$-projective modules.}
	\medskip
Recall that  $M$ is called $\ab$-projective provided $M$ is projective
in $\ab M$.  Clearly, this is equivalent to saying  that there 
exists an exact abelian subcategory $\Cal C$
of $\mod\Lambda$ which contains $M$ such that $M$ is projective inside $\Cal C.$ 
	\bigskip
{\bf Proposition 4.} {\it Let $M$ be an $\ab$-projective $\Lambda$-module 
and $\Gamma(M) = \End(M)^\op.$ Then }
$$
 \stat M = \cok M\quad \text{and}\quad \adstat M = \mod\Gamma(M).
$$
	\medskip
Proof. In general, $\stat M \subseteq \cok M.$ In order to show the equality, assume 
there is given an exact sequence
$$
 M'' \to M' \to N \to 0
$$
with $M', M'' \in \add M.$ The sequence shows that $N$ belongs to $\ab M,$ thus we deal with an
exact sequence in $\ab M.$ By assumption, $M$ is $\ab$-projective, thus the functor $\Hom(M,-)$
is right exact on $\ab M$. According to the implication (iv) $\implies$ (i) in Proposition 1
we see that $N$ belongs to $\stat M.$ This shows the first assertion.

In order to show the second assertion, let $X$ be in $\mod \Gamma(M)$. Let  
$$
 P_1 @>p>> P_0 @>q>> X @>>> 0
$$
be an exact sequence with finite-dimensional projective $\Gamma(M)$-modules $P_1, P_0$.
Since the the functor $M\otimes-$ is right exact, we obtain an exact sequence
$$
 M\otimes P_1 @>M\otimes p>>  M\otimes P_0 @>M\otimes q>>  M\otimes X @>>>  0.
$$
The two modules $M\otimes P_i$ belong to $\add M$, therefore the exact
sequence lies in $\ab M.$ 
Let us apply the functor $\Hom(M,-)$. Since $M$ is $\ab$-projective, the functor $\Hom(M,-)$
sends exact sequences in $\ab M$ to exact sequences, thus the sequence
$$
 \Hom(M,M\otimes P_1) @>>> \Hom(M,M\otimes P_0) @>>> \Hom(M, M\otimes X) @>>>  0.
$$
is exact. Altogether, there is the following commutative diagram with exact rows:
$$
\CD
 P_1 @>p>> P_0 @>q>> X @>>> 0 \cr
 @V\nu_{P_1}VV     @VV\nu_{P_0}V      @VV\nu_{X}V \cr
 \Hom(M,M\otimes P_1) @>>> \Hom(M,M\otimes P_0) @>>> \Hom(M, M\otimes X) @>>>  0.
\endCD
$$
The first two vertical maps are bijective, thus also $\nu_X$ is bijective.
\hfill$\square$
	\bigskip\bigskip
{\bf Proposition 5.} {\it Let $\Cal C$ be an exact abelian subcategory of $\mod\Lambda$,
let $C$ be a progenerator of $\Cal C$. Then
$C$ is an $\ab$-projective module and }
$$
 \Cal C = \stat C = \ab C.
$$
	\medskip 
Proof: First, let us show that $\ab C = \Cal C.$ Since $\Cal C$ is an exact
abelian subcategory which contains $C$, we have $\ab C \subseteq \Cal C.$
On the other hand, let $N$ be a module in $\Cal C$. Since $C$ is a progenerator in
$\Cal C$, there is an exact sequence $C'' \to C' \to N \to 0$ with $C', C''$
in $\add C$, thus $N$ is the cokernel of the map $C'' \to C'$ in $\add C$.
Since $\add C \subseteq \ab C$ and $\ab C$ is closed under cokernels, we see that
$N$ belongs to $\ab C.$ 

By assumption, $C$ is projective in $\Cal C$, thus $C$ is $\ab$-projective. \hfill $\square$
	\bigskip\bigskip 
{\bf 5. Triple modules.}
	\medskip
We consider indecomposable $\Lambda$-modules $M$ such that $\Gamma(M)$ is a 
Nakayama algebra of length 2. If $f$ is a non-zero nilpotent endomorphism of $M$, then
$f^2 = 0$, thus $\Im(f)\subseteq \Ker(f)$ and these submodules $\Im(f), \Ker(f)$ are
uniquely determined. Thus, $M$ has a uniquely determined filtration
$$
 0 = M_0 \subset M_1 \subseteq M_2 \subset M_3 = M,
$$
such that $M/M_2$ is isomorphic to $M_1$ (such an isomorphism is provided by $f$).
The module $M_1$ has to be a brick, since otherwise we obtain further endomorphisms of $M$.
We say that $M$ is a {\it triple module} provided also $M_2/M_1$ is isomorphic to $M_1.$
	\medskip
{\bf Proposition 6.} {\it If $M$ is a triple module, then $\stat M = \add M$ and this is
not an abelian category.}
	\medskip
Proof. Since $\Gamma(M)$ is a Nakayama algebra of length $2$, there are precisely
2 indecomposable $\Gamma(M)$-modules, thus, there are at most two isomorphism classes
of indecomposable $M$-static modules. Assume that there is an indecomposable $M$-static
module $N$ which is not isomorphic to $M$. According to Proposition 2, we can assume that $N = M/U$,
where $U$ is the image of an endomorphism of $M$ and such that the canonical projection
$M \to M/U$ is a right $M$-approximation. Since $M$ is a triple module, it follows that $U = M_1.$

Now $N = M/M_1$ has the submodule $M_2/M_1$. Since $M/M_2$ and $M_2/M_1$ are isomorphic, there is
a homomorphism $f\:M \to N = M/M_1$ with image $M_2/M_1$. Since the projection $p\:M \to M/M_1$ 
is a right $M$-approximation, there is a map $f'\;M \to M$ with $f = pf'.$ 
Since $\Gamma(M)$ is of length 2, either $f'$ is an automorphism or else $(f')^2 = 0.$ 
But $f'$ cannot be an automorphism, since the kernel of $p$ is $M_1$, whereas the kernel of $f$
is $M_2$. Also, $(f')^2 = 0$ is impossible, since in this case the image of $f'$ would be $M_1$
and then $pf' = 0,$ whereas $f\neq 0.$ This contradiction shows that the only indecomposable
$M$-static module is $M$. 

It remains to show that $\add M$ is not an abelian subcategory. Namely, consider a non-zero
nilpotent endomorphism $f$ of $M$. If $\add M$ is abelian, then it has to be a length category,
in particular, there has to exist a simple object in $\add M$, thus a brick. But $\add M$
has a unique indecomposable object, namely $M$, and by assumption, $M$ is not a brick. 
\hfill$\square$
	\bigskip\bigskip
  
{\bf 6. Nakayama algebras.}
	\medskip
{\bf Proposition 7.} 
{\it Any indecomposable module of a Nakayama algebra is $\ab$-projective.}
	\medskip 
Proof. Let $\Lambda$ be a Nakayama algebra  and $M$ an indecomposable $\Lambda$-module say of length $t$.
Let $J$ be the radical of $\Lambda$ and $\Lambda' = \Lambda/J^t$. Since $J^tM = 0$, we see that $M$
is a $\Lambda'$-module. The $\Lambda'$-modules form an exact abelian subcategory and as we have
seen, $M$ belongs to $\mod\Lambda'$, thus $\ab M \subseteq \mod\Lambda'$. Let 
$P = P_{\Lambda'}(M)$ be a
projective cover of $M$ considered as a $\Lambda'$-module. Then $P$ is an indecomposable $\Lambda'$-module.
Since any indecomposable $\Lambda'$-module has length at most $t$, the module $P$ has length at most $t$.
But $M$ is a factor module of $P$ and has length $t$. This shows that $M = P$ is projective in 
$\mod\Lambda'$ and therefore in $\ab M$. \hfill $\square$
	\bigskip 
Let us provide more details about the categories $\stat M$ and $\ab M$ for $M$ an indecomposable
$\Lambda$-module of length $t$, and $\Lambda$ a Nakayama algebra. Let us assume that the number of simple
$\Lambda$-modules is $s$. We may assume that $t \ge s$ (note that for $t\le s$, the module $M$ is
a brick, thus $\add M = \stat M = \ab M$). Since $\Lambda$ is a Nakayama algebra, the indecomposable
$\Lambda$-modules are uniquely determined by the length and the top (this is the isomorphism class
of a simple $\Lambda$-module). If $S$ is a simple module, we denote by $[i]S$ the indecomposable
$\Lambda$-module of length $i$ with top $S$. 
	\medskip
{\bf Proposition 8.} {\it Let $\Lambda$ be a Nakayama algebra with $s$ simple modules. Let $M$
be an indecomposable $\Lambda$-module of length $t \ge s$ with top $S$.
Let $e = \lceil {t \over s} \rceil$. 
Let $\Gamma(M) = \End(M)^\op.$ Then $\Gamma(M)$ is a local Nakayama algebra of length $e$.

{\rm(a)} 
If $s|t,$ then $\stat M = \ab M$ is equivalent to $\mod \Gamma(M)$ and 
the indecomposable modules in $\ab M$ are the modules
of the form $[i]S$ where $1\le i \le e.$ 

If $s > t$ and $t = (e-1)s + s_1$ with $1\le s_1 < s$, then $\ab M$ is equivalent to the module
category of a Nakayama algebra with two simple modules;
the modules in $\ab M$ have a filtration with factors of the form $[s_1]S$ and $[s]S/[s_1]S,$  and
an indecomposable module in $\ab M$ 
belongs to $\stat M$ if and only if it is either isomorphic to $M$ or of the form $[is]S$ with $1\le i < e.$ 

{\rm(b)} 
Always, $\ab M$ is equivalent to the module category of a Nakayama algebra $\Gamma'$ with at most two simple
modules and a unique indecomposable module which is both projective and injective
(namely the $\Gamma'$-module $\eta(M)$, where $\eta\:\ab M \to \mod \Gamma'$ is an equivalence). 

{\rm(c)} 
The category $\stat M$ is always an abelian subcategory, the embedding $\stat M \to \mod\Lambda$
is right exact, but usually not left exact. This embedding is exact if and only if $s|t.$
}
	\medskip 
Proof. First, assume that $s|t,$ thus $t = es.$ Consider the module $B = [s]S$, this is a brick and
$M$ has a filtration whose factors are all of the form $B$. Let $\Cal C$ be the subcategory of $\mod\Lambda$
whose objects are direct sums of modules of the form $[i]S$ with $1\le i \le e$. Then this is
an exact abelian subcategory, and $M$ belongs to $\Cal C$, thus $\ab M \subseteq \Cal C$. On the other hand,
every indecomposable module in $\Cal C$ is the cokernel of an endomorphism of $M$, thus 
$\Cal C \subseteq \cok M \subseteq \ab M$. It follows that $\ab M = \Cal C$. 
According to Proposition 7, $M$ is $\ab$-projective and according to Proposition 4, $\stat M = \cok M$.
and $\adstat M = \mod\Gamma(M)$. It follows that $\ab M$ is equivalent to the category 
$\mod\Gamma(M)$, the algebra $\Gamma(M)$ is a local Nakayama algebra, thus it has a unique
simple module, and the equivalence $\eta = \Hom(M,-)\:\ab M \to \mod\Gamma(M)$ sends $M$ to
$\eta(M)$, the only indecomposable projective (and also injective) $\Gamma(M)$-module.
	\medskip
Second, let $s > t$ and $t = (e-1)s + s_1$ for some $1\le s_1 < s;$ note that $e \ge 2$.
We consider $B_1 = [s_1]S$, this is a brick of length $s_1,$ and 
$B_2 = [s]S/[s_1]S$, this is a brick of length $s_2 = s-s_1.$
We consider the subcategory $\Cal F = \Cal F(B_1,B_2)$ of all $\Lambda$-modules with a filtration
with factors of the form $B_1,B_2.$ This is an exact abelian subcategory.
The module $M$ has a filtration with $e$ factors $B_1$ and $e-1$ factors $B_2$, thus it
belongs to $\Cal F$; and therefore $\ab M \subseteq \Cal F.$ 
Of course, $\Cal F$ is equivalent to the module category of a
Nakayama algebra with precisely two simple modules, and we denote by $\Cal C$ its subcategory 
of all direct sums of indecomposable objects which have a filtration with factors $B_1, B_2$, such 
that there are at most $e-1$ factors of the form $B_2.$ We have $\ab M \subseteq \Cal C.$ 

Let us show that $\ab M = \Cal C$. The modules $N$ in $\Cal C$ with $\Hom(B_1,N) \neq 0 \neq \Hom(N,B_1)$
are images of endomorphisms of $M$, those with $\Hom(B_1,N) \neq 0 \neq \Hom(N,B_2)$ are
kernels of endomorphisms of $M$, those with $\Hom(B_2,N) \neq 0 \neq \Hom(N,B_1)$ are
cokernels of endomorphisms of $M$. Finally, the modules $N$ in $\Cal C$ with $\Hom(B_2,N) \neq 0 
\neq \Hom(N,B_2)$ are obtained as kernels of maps $M/B_1 \to M$. 

Let us denote by $M'$ the kernel of a non-zero map $M \to B_1$. Then $M\oplus M'$ is a progenerator 
of $\Cal C$. If $\Gamma' = \End(M\oplus M')^\op,$ then $\Cal C$ is equivalent to $\mod\Gamma'$
under the functor $\eta = \Hom(M\oplus M',-):\Cal C \to \mod \Gamma'.$ The algebra
$\Gamma'$ is a Nakayama algebra with precisely two simple modules and $\eta(M)$ is indecomposable
and both projective and injective as a $\Gamma'$-module. 

The modules in $\stat M$ are precisely the cokernels of endomorphisms of $M$, see Proposition 4.
The indecomposable modules $N$ which belong to $\stat M$ and are not isomorphic to $M$
are precisely the cokernels of non-zero endomorphisms of $M$, 
thus the modules $N$ in $\Cal C$ with $\Hom(B_2,N) \neq 0 \neq \Hom(N,B_1)$, or, equivalently,
those of the form $N = [is]S$ for some $1\le i \le e-1$.

Finally, always $\stat M$ is equivalent to the abelian category $\adstat M = \mod\Gamma(M)$,
thus it is an abelian subcategory of $\mod\Lambda$. The equivalence is given by the 
functor $M\otimes - \:\mod\Gamma(M) \to \mod\Lambda$, this is a right exact (but usually not
left exact) functor. In case $s|t$, the equality $\stat M = \ab M$ shows that $\stat M$
is an exact abelian subcategory. If $s$ does not divide $t$, say 
let $t = (e-1)s + s_1$, there is a non-zero map $f\:[s]S \to M$, its kernel in $\mod\Lambda$
is $B_2$, thus does not belong to $\stat M$. This shows that $\stat M$ is not closed under kernels.
\hfill $\square$
	\bigskip
Let us exhibit an  example. We consider the path algebra of the quiver
$$
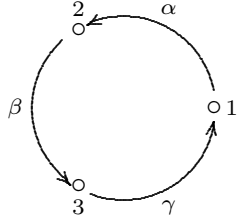
\hbox{\beginpicture
	\setcoordinatesystem units <1cm,1cm>
\multiput{} at -1 -1  1 1 /
\multiput{$\circ$} at  1 0  -.5 .866  -.5 -.866 /
\circulararc 100 degrees from 1 0.2 center at 0 0
\circulararc -100 degrees from 1 -.2 center at 0 0
\circulararc 100 degrees from -.68 0.75 center at 0 0
\arr{0.99 -.2}{1 -.15}
\arr{-.38 .93}{-.4 .92}
\arr{-.64 -.8}{-.59 -.85}
\put{$\ssize 1$} at 1.2 0
\put{$\ssize 2$} at -.5 1.1
\put{$\ssize 3$} at -.5 -1.1
\put{$\ssize \alpha$} at 0.5 1.1
\put{$\ssize \beta$} at -1.2 0 
\put{$\ssize \gamma$} at 0.5 -1.1
\put{with relations} [l] at  2 0.3
\put{$\beta\alpha\gamma\beta\alpha\gamma\beta\alpha = 0 = \gamma\beta\alpha\gamma\beta\alpha\gamma.$} 
  [l] at 2 -.2
\endpicture} 
$$
It is a Nakayama algebra $\Lambda$ with Kupisch series $(8,8,7)$ (the Kupisch
series of a Nakayama algebra records the numbers $(p_1,\dots,p_s)$, where $p_i$ is the
length of the projective cover $P(i)$ of the simple module $S(i)$). 
The following picture shows the Auslander-Reiten quiver of $\Lambda$. We choose $M = P(1)$.
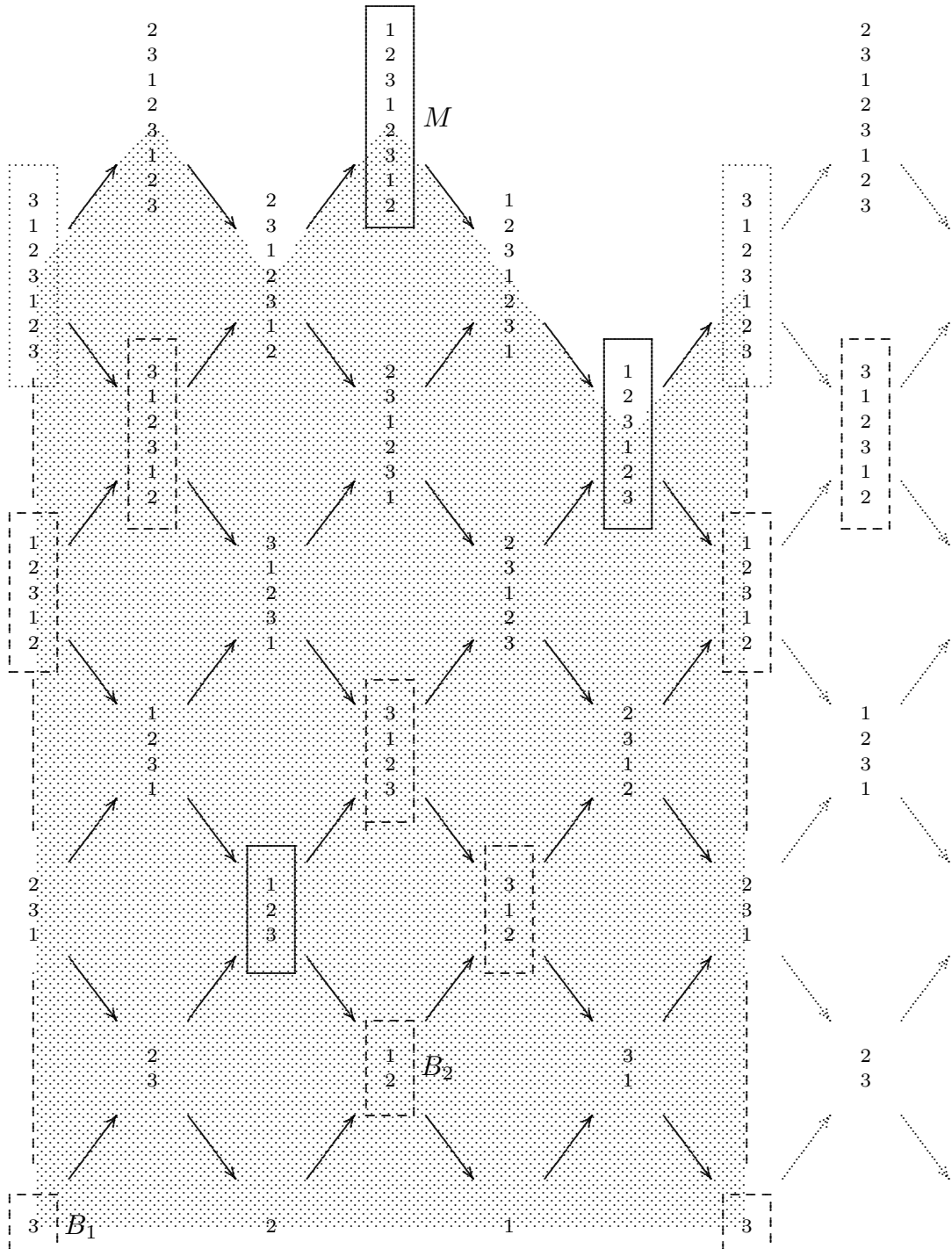
\vfill\eject 
$$
\hbox{\beginpicture
	\setcoordinatesystem units <1.5cm,2cm>
\plot 2.8 6.3  3.2 6.3  3.2 7.7  2.8 7.7  2.8 6.3 /
\plot 4.8 4.4  5.2 4.4  5.2 5.6  4.8 5.6  4.8 4.4 /
\plot 1.8 1.6  2.2 1.6  2.2 2.4  1.8 2.4  1.8 1.6 /

\setdashes <1mm>
\plot 5.8 -.2  6.2 -.2  6.2 .2  5.8 .2  5.8 -.2 /
\plot -.2 -.2  0.2 -.2  0.2 .2  -.2 .2  -.2 -.2 /

\plot 2.8 .7  3.2 .7  3.2 1.3  2.8 1.3  2.8 .7 /

\plot 3.8 1.6  4.2 1.6  4.2 2.4  3.8 2.4  3.8 1.6 /

\plot 2.8 2.55  3.2 2.55  3.2 3.45  2.8 3.45  2.8 2.5 /

\plot 5.8 3.5  6.2 3.5  6.2 4.5  5.8 4.5  5.8 3.5 /
\plot -.2 3.5  0.2 3.5  0.2 4.5  -.2 4.5  -.2 3.5 /

\plot .8 4.4  1.2 4.4  1.2 5.6  .8 5.6  .8 4.4 /
\plot 6.8 4.4  7.2 4.4  7.2 5.6  6.8 5.6  6.8 4.4 /

\setdots <.8mm>
\plot 5.8 5.3  6.2 5.3  6.2 6.7  5.8 6.7  5.8 5.3 /
\plot -.2 5.3  0.2 5.3  0.2 6.7  -.2 6.7  -.2 5.3 /

\setsolid
\put{$M$} at 3.4 7
\multiput{$\smallmatrix 2\cr 3\cr1\cr 2\cr 3\cr1\cr 2\cr 3\endsmallmatrix$} at 1 7  7 7 /
\multiput{$\smallmatrix 1\cr 2\cr3\cr 1\cr 2\cr3\cr 1\cr 2\endsmallmatrix$} at 3 7 /

\multiput{$\smallmatrix 3\cr1\cr 2\cr 3\cr1\cr 2\cr 3\endsmallmatrix$} at 0 6  6 6 /
\multiput{$\smallmatrix 2\cr3\cr 1\cr 2\cr3\cr 1\cr 2\endsmallmatrix$} at 2 6 /
\multiput{$\smallmatrix 1\cr2\cr 3\cr 1\cr2\cr 3\cr 1\endsmallmatrix$} at 4 6 /

\multiput{$\smallmatrix  3\cr1\cr 2\cr 3\cr1\cr 2\endsmallmatrix$} at 1 5  7 5 /
\multiput{$\smallmatrix  2\cr3\cr 1\cr 2\cr3\cr 1\endsmallmatrix$} at 3 5 /
\multiput{$\smallmatrix  1\cr2\cr 3\cr 1\cr2\cr 3\endsmallmatrix$} at 5 5 /

\multiput{$\smallmatrix 1\cr 2\cr 3\cr1\cr2\endsmallmatrix$} at 0 4  6 4 /
\multiput{$\smallmatrix 3\cr 1\cr 2\cr3\cr1\endsmallmatrix$} at 2 4 /
\multiput{$\smallmatrix 2\cr 3\cr 1\cr2\cr3\endsmallmatrix$} at 4 4 /

\multiput{$\smallmatrix 1\cr 2\cr 3\cr1\endsmallmatrix$} at 1 3  7 3 /
\multiput{$\smallmatrix 3\cr 1\cr 2\cr3\endsmallmatrix$} at 3 3 /
\multiput{$\smallmatrix 2\cr 3\cr 1\cr2\endsmallmatrix$} at 5 3 /

\multiput{$\smallmatrix 2\cr 3\cr1 \endsmallmatrix$} at 0 2  6 2 /
\multiput{$\smallmatrix 1\cr 2\cr3 \endsmallmatrix$} at 2 2 /
\multiput{$\smallmatrix 3\cr 1\cr2 \endsmallmatrix$} at 4 2 /

\multiput{$\smallmatrix 2\cr 3\endsmallmatrix$} at 1 1  7 1 /
\multiput{$\smallmatrix 1\cr 2\endsmallmatrix$} at 3 1 /
\multiput{$\smallmatrix 3\cr 1\endsmallmatrix$} at 5 1 /

\multiput{$\smallmatrix 3\endsmallmatrix$} at 0 0  6 0 /
\multiput{$\smallmatrix 2\endsmallmatrix$} at 2 0 /
\multiput{$\smallmatrix 1\endsmallmatrix$} at 4 0 /

\arr{0.3 0.3}{0.7 0.7}
\arr{1.3 0.7}{1.7 0.3}
\arr{2.3 0.3}{2.7 0.7}
\arr{3.3 0.7}{3.7 0.3}
\arr{4.3 0.3}{4.7 0.7}
\arr{5.3 0.7}{5.7 0.3}

\arr{0.3 1.7}{0.7 1.3}
\arr{1.3 1.3}{1.7 1.7}
\arr{2.3 1.7}{2.7 1.3}
\arr{3.3 1.3}{3.7 1.7}
\arr{4.3 1.7}{4.7 1.3}
\arr{5.3 1.3}{5.7 1.7}

\arr{0.3 2.3}{0.7 2.7}
\arr{1.3 2.7}{1.7 2.3}
\arr{2.3 2.3}{2.7 2.7}
\arr{3.3 2.7}{3.7 2.3}
\arr{4.3 2.3}{4.7 2.7}
\arr{5.3 2.7}{5.7 2.3}

\arr{0.3 3.7}{0.7 3.3}
\arr{1.3 3.3}{1.7 3.7}
\arr{2.3 3.7}{2.7 3.3}
\arr{3.3 3.3}{3.7 3.7}
\arr{4.3 3.7}{4.7 3.3}
\arr{5.3 3.3}{5.7 3.7}

\arr{0.3 4.3}{0.7 4.7}
\arr{1.3 4.7}{1.7 4.3}
\arr{2.3 4.3}{2.7 4.7}
\arr{3.3 4.7}{3.7 4.3}
\arr{4.3 4.3}{4.7 4.7}
\arr{5.3 4.7}{5.7 4.3}

\arr{0.3 5.7}{0.7 5.3}
\arr{1.3 5.3}{1.7 5.7}
\arr{2.3 5.7}{2.7 5.3}
\arr{3.3 5.3}{3.7 5.7}
\arr{4.3 5.7}{4.7 5.3}
\arr{5.3 5.3}{5.7 5.7}

\arr{0.3 6.3}{0.7 6.7}
\arr{1.3 6.7}{1.7 6.3}
\arr{2.3 6.3}{2.7 6.7}
\arr{3.3 6.7}{3.7 6.3}

\setdots <.5mm>
\arr{6.3 0.3}{6.7 0.7}
\arr{7.3 0.7}{7.7 0.3}
\arr{6.3 1.7}{6.7 1.3}
\arr{7.3 1.3}{7.7 1.7}
\arr{6.3 2.3}{6.7 2.7}
\arr{7.3 2.7}{7.7 2.3}
\arr{6.3 3.7}{6.7 3.3}
\arr{7.3 3.3}{7.7 3.7}

\arr{6.3 4.3}{6.7 4.7}
\arr{7.3 4.7}{7.7 4.3}

\arr{6.3 5.7}{6.7 5.3}
\arr{7.3 5.3}{7.7 5.7}

\arr{6.3 6.3}{6.7 6.7}
\arr{7.3 6.7}{7.7 6.3}

\setdashes <1mm> 
\plot 0 0.4  0 1.6 /
\plot 0 2.5  0 3.5 /
\plot 0 4.6  0 5.4 /
\plot 6 0.4  6 1.6 /
\plot 6 2.5  6 3.5 /
\plot 6 4.6  6 5.4 /

\setshadegrid span <1.5pt>
\hshade  0 0 6 <,,z,z> 5  0 6 <,,z,z> 5.01 0 5 <,,z,z> 6 0  4 
  <,,z,z> 6.01 0 2 <,,z,z> 7 1  1 /
\hshade   6.0 2 4 <,,z,z> 7 3 3 /
\hshade   5.0 5 6 <,,z,z> 6 6 6 /
\put{$B_1$} at .4 0
\put{$B_2$} at 3.4 1

\endpicture}
$$
Let $B_1 = S(1)$ and $B_2$ be the indecomposable module of length
$2$ with socle $S(2)$ and top $S(1)$. Then
$$
  \stat M \subset \ab M \subset \Cal F(B_1,B_2).
$$
The modules in $\stat M$ are marked by a solid frame, those
belonging to $\ab M$, but not to $\stat M$ have a dashed frame. 
There is one additional indecomposable 
module which belongs to $\Cal F(B_1,B_2)$, but not to $\ab M$, 
it has a dotted frame.
   \bigskip\bigskip 
{\bf 7. Proof of Theorem 1 and Theorem 2.}
   \medskip 
Let $\Lambda$ be a finite-dimensional  hereditary $k$-algebra. 
We may assume in addition that $\Lambda$ is connected (this means that $0$ and $1$ are
the only central idempotents). 
Let us recall relevant properties of $\mod\Lambda$. 

Let $K_0(\Lambda)$ be the Grothendieck group of $\mod\Lambda$, this is the free abelian
group with basis the isomorphism classes of all the finite-dimensional $\Lambda$-modules
and with relations of the form $[X]-[Y]+[Z]$ provided there exists an exact
sequence of the form $0 \to X \to Y \to Z \to 0$. Given a $\Lambda$-module $X$, 
the residue class
of its isomorphism class is denoted by $\bdim X$. Note that $K_0(\Lambda)$ is
the free abelian group with basis the elements $\bdim S$, where $S$ runs through the simple
$\Lambda$-modules. Since $\Lambda$ is hereditary, there is a (uniquely determined) bilinear
form $\langle-,-\rangle$ on $K_0(\Lambda)$ such that
$$
  \langle \bdim X, \bdim Y\rangle = \dim\Hom(X,Y) - \dim\Ext^1(X,Y).
$$
The corresponding quadratic form $q_\Lambda$ defined by
$q_\Lambda(x) = \langle x,y\rangle$ for $x\in K_0(\Lambda)$  
it is called the Euler form for $\Lambda$.

The algebra $\Lambda$ is said to be {\it tame} provided the form $q_\Lambda$ is
positive semidefinite, otherwise {\it wild}. One may attach to $\Lambda$ a 
{\it valued quiver} $Q(\Lambda)$ 
as follows: it has as vertices the isomorphism classes
$[S]$ of the simple $\Lambda$-modules and there is an arrow $[S] \to [S']$ provided $\Ext(S,S') \neq 0$.
Since we assume that $\Lambda$ is finite-dimensional, the quiver $Q(\Lambda)$ is finite and 
directed (the latter means that the simple modules can be labeled 
$S(i)$ such that the existence of an arrow $[S(i)]\to [S(j)]$ implies that $i > j$). 
We endow $Q(\Lambda)$ with a valuation as 
follows: Given an arrow $[S] \to [S']$, consider $\Ext(S,S')$ as a left $\End(S)^\op$-module
and also as a left $\End(S')$-module and put
$$
 v([S],[S']) = \dim {}_{\End(S)^\op}\Ext(S,S')\times\dim {}_{\End(S')}\Ext(S,S')
$$
provided $v([S],[S']) > 1$. The algebra $\Lambda$ is tame if and only if the underlying valued graph
(obtained by forgetting the orientation of the edges) is either a Dynkin diagram $\Bbb A_n, \Bbb B_n,\dots,
\Bbb G_2$ or a Euclidean diagram $\widetilde{\Bbb A}_n, \dots, \widetilde {\Bbb G}_{21},
\widetilde {\Bbb G}_{22}.$  

If $\Lambda$ is tame, the structure of
the category is completely known:  the category
$\mod\Lambda$ consists of two directed components and a family of separating tubes. This
implies: If $M$ is any
indecomposable $\Lambda$-module, then $\Gamma(M)$ is a (local) Nakayama algebra
say of length $e(M).$ 

If $e(M) = 1$, then $\add M$ is an abelian subcategory, thus
$\add M = \ab M$, and this is a semisimple abelian category. In particular,
$M$ is $\ab$-projective.  

Assume now that $e(M) \ge 2$. Then $M$ belongs to an Auslander-Reiten
component $\Bbb C$ 
which is a tube (i.e\. of the form $\Bbb Z\Bbb A_\infty/n$ for some $n$). 
The modules which are direct sums of indecomposables in $\Bbb C$ form an exact
abelian subcategory which we also denote by $\Cal C$. 
For any natural number $b$, the subcategory $\Cal C_b$
of $\Cal C$ given by all objects
in $\Cal C$ with (relative) Loewy length at most $b$ is an exact abelian subcategory
which is equivalent to the module category of a Nakayama algebra. This shows that
$\ab M$ can be considered as an exact abelian subcategory of the module
category of a Nakayama algebra. It follows from Proposition 7 that $M$ is $\ab$-projective.  

This completes the proof of the implication (i) $\implies$ (ii) in Theorem 1.
The implication (ii) $\implies$ (iii) has been shown in Proposition 5, the
implication (iii) $\implies$ (iv) is trivial since $\stat M$ is equivalent
as a category to $\adstat M$. According to Proposition 4, we have $\stat M = \cok M,$
thus condition (v) is satisfied.

Now assume that $\Lambda$ is wild. According to [R], $\Lambda$ is strictly wild (a $k$-algebra is
{\it strictly wild} provided  there is
a finite extension field $k'$ of $k$ such that for any finite-dimensional $k'$-algebra
$\Gamma$, there is a full and exact embedding $\mod \Gamma \to \mod\Lambda$).

It remains to show the implications (iv) $\implies$ (i) and 
(v) $\implies$ (i) in Theorem 1, and the implications (ii) $\implies$ (i)
and (iii) $\implies$ (i) in Theorem 2. 

Assume that $\Lambda$ satisfies one of the
conditions (iv) or (v) of Theorem 1. We claim that $\Lambda$ cannot be wild.
Namely, if $\Lambda$ is wild, then according to the implication (i) $\implies$ (ii) of
Theorem 2, there exists an indecomposable module $M$ with $\Gamma(M)$
not a division ring, such that $\stat M = \add M$. But $\add M$ cannot be an
abelian category, this contradicts (iv). 
Also, since $\Gamma(M)$ is not a division ring,
$\add M$ is a proper subcategory of $\cok M$, a contradiction to (v). 

Assume now there exists an indecomposable $\Lambda$-module $M$ such that $M$ is not a division ring
and $\stat M = \add M.$ We claim that $\Lambda$ has to be wild. Namely, if $\Lambda$ is tame, then 
the implication (i) $\implies$ (iv) of Theorem 1
asserts that $\stat M$ is abelian. But if $M$ is indecomposable and not a brick, then $\add M$ is not
abelian. Similarly, assume that there is a finite extension field $k'$ of $k$ such that for
any finite-dimensional $k$-algebra $\Gamma$, there is a $\Lambda$-module $M$ such that $\stat M$ is
equivalent to $\mod\Gamma.$ In particular, take any representation-finite $k'$-algebra $\Gamma$
and a $\Lambda$-module $M$ such that $\stat M$ is
equivalent to $\mod\Gamma.$ Then $\Lambda$ has to be wild. Namely, if $\Lambda$ would be tame,
then the implication (i) $\implies$ (iii) of Theorem 1 asserts that $\stat M$ is equivalent to the
module category of a Nakayama algebra, in particular, that $\stat M$ has only finitely many isomorphism
classes of indecomposable objects. 
\hfill$\square$
	\bigskip
{\bf Remark.} Assume again that $\Lambda$ is a tame hereditary algebra. We have seen
that for $M$ indecomposable, there is a full embedding of $\mod\Gamma(M)$ into $\mod\Lambda$.
If $M$ is not indecomposable, then usually $\mod\Gamma(M)$ will not be equivalent to a 
subcategory of $\mod\Lambda$.
	\smallskip 
{\bf Example.} 
Consider $\Bbb D_4$ with subspace orientation and let 
$M = S(0)\oplus \tau^{-1}S(0)$,where $\tau$ is the Auslander-Reiten translation. 
Then $\Gamma(M)$ is the Kronecker algebra, 
in particular, representation-infinite. This shows that $\mod\Gamma(M)$ cannot be
a  subcategory of $\mod\Lambda$. 

Let us exhibit a $\Lambda$-module $N$ which is not $M$-static. 
Take a map $S(0) \to \tau^{-1}S(0)$ in general position, thus the cokernel $N$ is the
injective hull of $S(0)$. But the minimal right $M$-approximation of $N$ is of the
form $M^2\to N$ with kernel $\Omega_M(N)$ being 
the direct sum of the indecomposable projective modules
of length 2. Of course, $\Omega_M(N)$ is not generated by $M$. Thus, $N$ is not $M$-static.
   \bigskip\bigskip

{\bf Acknowledgement.} 
This work is funded by the Deanship of Scientific Research, 
King Abdulaziz University, under grant No. 2-130/1434/HiCi. 
The authors, therefore, acknowledge technical and financial support of KAU.
    \bigskip\bigskip  
\vfill\eject 
{\bf References.}
\frenchspacing
     \medskip 
\item{[A]}J.~L.~Alperin: 
    Static modules and non-normal Clifford theory, J. Aust. Math. Soc. (Series A) 49, (1990), 347--353.
\item{[ARS]} M.~Auslander, I.~Reiten, S.~Smal\o: Representation Theory of  
    Artin Algebras. Cambridge Studies in Advanced Mathematics 36. Cambridge University 
    Press. 1997.
\item{[DR1]} V.~Dlab, C.~M.~Ringel:
    Indecomposable representations of graphs and algebras. Mem. Amer. Math. Soc. 173 (1976). 
\item{[DR2]}V.~Dlab, C.~M.~Ringel:  The global dimension of the endomorphism ring of a 
    generator-cogenerator for a hereditary artin algebra.
    Mathematical Reports of the Academy of Science of the Royal Society of Canada. Vol 30 (2008), Nr.3, 89--96. 
\item{[N1]}S.~K.~Nauman, Static modules and stable Clifford theory, J. Algebra,128, (1990), 497--509.
\item{[N2]} S.~K.~Nauman, Static modules, Morita contexts, and equivalences. J. Algebra 135 (1990), 192--202

\item{[R]} C.~M.~Ringel: Representations of K-species and bimodules. J. Algebra 41 (1976),    269-302.
\item{[W]}R.~Wisbauer: Static modules and equivalences, 
   in {\it Interactions between Ring Theory and Representation Theory}, V.~Oystaeyen, M.~Saorin (Ed.),
  Marcel Dekker, LNPAM 210 (2000), 423--449.
	   \bigskip\medskip 

{\rmk
\noindent 
Mustafa A. A. Obaid: E-mail: {\ttk drmobaid\@yahoo.com}

\noindent 
S. Khalid Nauman: E-mail: {\ttk snauman\@kau.edu.sa}

\noindent 
Wafaa  M. Fakieh: E-mail: {\ttk wafaa.fakieh\@hotmail.com}

\noindent 
Claus Michael Ringel: E-mail: {\ttk ringel\@math.uni-bielefeld.de}
	\smallskip 

King Abdulaziz University, Faculty of Science, \par 
P.O.Box 80203, Jeddah 21589, Saudi Arabia
}
\bye